\documentclass[11pt,twoside]{amsart}

\author{Valery Alexeev} 
\address{Department of Mathematics\\
University of Georgia\\
Athens, GA 30602}
\email{valery@math.uga.edu}

\renewcommand{\theenumi}{(\roman{enumi})}
\newcommand{\bP}{{\mathbb P}}
\newcommand{\bQ}{{\mathbb Q}}
\newcommand{\bZ}{{\mathbb Z}}

\newcommand{\bC}{{\mathbb C}}

\newcommand{\cL}{{\mathcal L}}

\newcommand{\cM}{{\mathcal M}}

\newcommand{\cO}{{\mathcal O}}

\newcommand{\wP}{\widetilde P}

\newcommand{\Supp}{\operatorname{Supp}}
\newcommand{\Proj}{\operatorname{Proj}}

\newcommand{\chr}{\operatorname{char}}

\theoremstyle{plain}
\newtheorem{thm}{Theorem}[section]

\newtheorem{lem}[thm]{Lemma}
\newtheorem{cor}[thm]{Corollary}

\newtheorem{conj}[thm]{Conjecture}

\theoremstyle{definition}
\newtheorem{defn}[thm]{Definition}

\newtheorem{saynum}[thm]{}
\newtheorem{exmp}[thm]{Example}

\newtheorem{rem}[thm]{Remark}

\newtheorem{ack}{Acknowledgments}   

\theoremstyle{remark}

\newenvironment{say}{}{}

\newcommand{\cond}{\operatorname{cond}}

\begin{document}
\bibliographystyle{amsalpha+}
\title[Stable pairs]%
{Log canonical singularities and complete moduli of stable pairs}
\date{August 15, 1996} 
\maketitle

{\bf
\begin{center}
  \begin{tabular}{@{}l@{}}
    Preliminary version
  \end{tabular}
\end{center}
}

\tableofcontents

\setcounter{section}{-1}

\section{Introduction}

\begin{saynum}
  This paper consists of two parts. In the first part, assuming the
  log Minimal Model Program (which is currently only known to be true
  in $\dim\le3$), we construct the complete moduli of ``stable pairs''
  $(X,B)$ of projective schemes with divisors that generalize the
  moduli space of $n$-pointed stable curves $M_{g,n}$ to arbitrary
  dimension. The construction itself is a direct generalization of
  that of \cite{Alexeev_Mgn} where it was given in the case of
  surfaces, and is based in part on ideas from
  \cite{KollarShepherdBarron88,Kollar_ProjCompleteModuli,Viehweg95}.
\end{saynum}

\begin{saynum}
  In the second part of the paper we study the singularities of
  stable quasiabelian varieties and stable quasiabelian pairs $(X,B)$
  that appear in \cite{AlexeevNakamura96} as limits of abelian
  varieties.  We show that the singularities are semi log canonical.
  This implies, via Koll\'ar's Ampleness Lemma, that over $\mathbb C$
  if there exists a compactification of the moduli space $A_g$ of
  principally polarized abelian varieties by stable quasiabelian
  pairs, then it is in fact projective.
  
  We give more examples of situations where log canonical singularities
  appear naturally in connection with complete moduli problems. One of
  them is the minimal and toroidal compactifications of quotients
  $D/\Gamma$ of bounded symmetric domains by arithmetic groups. We
  point out the fact, which could be obvious to specialists had they
  known the definitions, that they all have log canonical singularities
  and that the minimal (=Baily-Borel) compactification is the log
  canonical model of any toroidal compactification when the
  ``boundary'' divisor $B$ is correctly defined. 
\end{saynum}

\begin{ack}
  I would like to thank Professors V. Batyrev, I. Nakamura, Y.
  Kawamata, J.  Koll\'ar, S. Mori and R. Varley for very helpful
  conversations.
\end{ack}


\section{Definitions for singularities}
\label{sec:Definitions for singularities}

\begin{defn}
  Let $X$ be a normal variety (not necessarily irreducible) defined
  over an algebraically closed field $k$ of any characteristic, and
  let $B_1\dots B_m$ be distinct reduced divisors on $X$. Denote $\sum
  B_j$ by $B$. Let $i:U\hookrightarrow X$ be the inclusion of the
  nonsingular part and denote 
  \begin{displaymath}
    \cO(N(K+B))= i_*\cO(N(K_U+B\vert_U)),
  \end{displaymath}
  where $\cO(K_U)$ is the canonical sheaf, the top exterior power of
  $\Omega^1_U$, and $N$ is an integer.
  
  One says that the pair $(X,B)$ has {\em log canonical
    singularities\/} if 
  \begin{enumerate}
  \item $\cO(N(K+B))$ is invertible for some $N>0$ (one then says
    that $K+B$ is $\bQ$-Cartier).
  \item for any birational morphism from a normal scheme $f:Y\to X$
    one has
    \begin{displaymath}
      f_*\cO_Y\big( N(K_Y+f^{-1}B+\sum E_i) \big) 
      = \cO_X\big(N(K_X+B)\big),
    \end{displaymath}
    where $E_i$ are exceptional divisors of $f$.
  \end{enumerate}
\end{defn}

\begin{rem}
  The above definition can be formulated also for the case of a
  divisor $B=\sum b_jB_j$ with rational coefficients $b_j$ by
  requiring $N$ to be divisible enough.
\end{rem}
  
\begin{saynum}
  An equivalent way would be to use {\em log codiscrepancies\/} -- the
  coefficients $a_i$ appearing in the following natural formula:
  \begin{displaymath}
    f^*(K_X+B)=K_Y+f^{-1}B+\sum a_iE_i
  \end{displaymath}
  The log codiscrepancies depend only on the divisors $E_i$
  themselves, i.e. the corresponding discrete valuations of the
  function field, and not on the model $Y$ chosen. Indeed, every two
  models $Y_1$ and $Y_2$ are comparable since they are both dominated
  by a third normal variety $Y_3$ -- take for example the component of
  the normalization of $Y_1\underset{Y}{\times}Y_2$ which dominates
  $Y$.
\end{saynum}

\begin{defn}
  The singularities are log canonical if all log codiscrepancies are
  $\le1$. They are {\em log terminal\/} if $a_i<1$, {\em klt\/} if
  $a_i<1$ and $b_j<1$. And they are {\em canonical (resp. terminal)\/}
  if $B$ is empty and $a_i\le0$ (resp. $a_i<0$).
\end{defn}

\begin{rem}
  In the above definition one usually assumes $Y$ to be non-singular,
  and then one needs the embedded resolution of singularities and
  hence characteristic $0$. This does not appear to be necessary.
  Still, without the resolution of singularities the situation becomes
  somewhat cumbersome. For example, it is not absolutely obvious that
  the next definition is equivalent to, or even implies, the previous
  one (this {\em is} obvious with resolution of singularities).
\end{rem}

\begin{defn}
  Let $(X,B)$ be as above. We say that this pair has {\em pre log
    canonical singularities\/} if there exists a proper birational
  morphism from a nonsingular variety $f:Y\to X$ such that
  \begin{enumerate}
  \item $\cO(N(K+B))$ is invertible for some $N>0$.
  \item the exceptional set of $f$ is a union of divisors $E_i$.
  \item $\cup f^{-1}B_j\cup E_i$ has normal crossings.
  \item $f_*\cO_Y(N(K_Y+f^{-1}B+\sum E_i)) = \cO_X(N(K_X+B))$.
  \end{enumerate}
\end{defn}

\begin{say}
  Another important class of singularities is the following.
\end{say}

\begin{defn}\label{defn:semi_log_canonical}
  Let $X$ be a reduced variety (not necessarily irreducible) defined
  over an algebraically closed field $k$ of any characteristic, and
  let $B_1\dots B_m$ be distinct reduced divisors on $X$, denote $\sum
  B_j$ by $B$. In addition, assume that $X$ is quasi-projective over
  $k$. 
  
  Let $i:U\hookrightarrow X$ be the union of the open locus of
  Gorenstein points of $X$ not contained in $B$ and the nonsingular
  locus of $X$.  Denote
  \begin{displaymath}
    \cO(N(K+B))= i_*\cO(N(K_U+B\vert_U)),
  \end{displaymath}
  where $N$ is an integer, and $\cO(K_U)$ is the restriction of the
  dualizing sheaf $\omega_{\overline{U}}$ of a projective closure of
  $U$ (\cite[III.7]{Hartshorne77}).
  
  One says that the pair $(X,B)$ has {\em semi log canonical
    singularities\/} if
  \begin{enumerate}
  \item $X$ satisfies the Serre's condition $S_2$.
  \item $X$ is Gorenstein in codimension $1$.
  \item none of the irreducible components of $B_j$ is contained in
  the singular locus of $X$.
  \item the closed subscheme $\cond(\nu)$ of $X^{\nu}$ corresponding
    to the conductor of normalization $\nu:X^{\nu}\to X$ is a union of
    reduced divisors.
  \item $\cO(N(K+B))$ is invertible for some $N>0$.
  \item the pair $(X^{\nu},\nu^{-1}B+\cond(\nu))$ has log canonical
    singularities.
  \end{enumerate}
\end{defn}

\begin{say}
  In the same way as above, one can define pre semi log canonical
  singularities.
\end{say}

\begin{rem}
  We note a certain lack of symmetry in the definitions of
  $\cO(N(K+B))$ for log and semi log canonical cases. However, they
  coincide if $X$ is both normal and quasi-projective.
\end{rem}

\begin{defn}
  Under the assumptions above, we will say that $K+B$ is ample if the
  sheaf $\cO_X(N(K+B))$, equivalently
  $\cO_{X^{\nu}}(N(K+B+\cond(\nu)))$, is an ample invertible sheaf for
  some integer $N>0$. In this case the pair $(X,B)$ is called the {\em
    log canonical model}.
\end{defn}

\begin{say}
  Let us now try to see what is the most general situation where the
  previous definitions still work. The main thing to understand is the
  canonical sheaf. The rest transfers over in a pretty straightforward
  way.
\end{say}
\smallskip

\begin{say}
  Let us fix a regular Noetherian scheme $S$ (for example spectrum of
  $\mathbb Z$ or a DVR) and consider a reduced scheme $X$ flat and of
  finite type over $S$.  Let us assume that $\pi:X\to S$ is smooth in
  codimension 1.  Denoting by $i:U\to X$ the embedding of this smooth
  locus, we can set
  \begin{displaymath}
    \cO_X(K_{X/S})=i_* \cO_X(K_{U/S}),
  \end{displaymath}
  where $\cO_X(K_{U/S})$ is the top exterior power of
  $\Omega^1_{U/S}$.  We can now define (pre) log canonical
  singularities of a pair $(X,B=\sum B_j)$, where $B_j\subset X$ are
  closed codimension 1 subschemes of $X$, by copying the definitions
  from section \ref{sec:Moduli of stable pairs in general}. In
  particular, for log canonical singularities we require $X$ to be
  normal.

  For (pre) semi log canonical we need to assume that $X/S$ is
  quasi-projective and that the normalization $X^{\nu}/S$ is smooth in
  codimension 1.
  
  As one can see, in these definitions we use the regular scheme $S$
  only as ``the beginning of coordinates'', something to start
  measuring from.
\end{say}

\medskip

\begin{say}
  Let us push the limits even little further. Clearly, the definition
  of (pre) log canonical singularities is stable under \'etale maps.
  Therefore, they transfer directly to algebraic spaces and algebraic
  stacks. If $R_X\overset{\to}{\to} U_X$ is an equivalence relation or
  a groupoid defining $X$, and $R_{B_j},U_{B_j}$ are the closed
  subschemes corresponding to $B_j$, then we say that the pair $(X,B)$
  has (pre) log canonical singularities if the same holds for
  $(U_X,U_B)$.
\end{say}


\section{Moduli of stable pairs in general}
\label{sec:Moduli of stable pairs in general}

\begin{say}
  The purpose of this section is to describe a construction of
  complete and projective moduli spaces for stable $n$-dimensional
  pairs which generalize the usual moduli of stable $n$-pointed
  curves.  This will be done assuming a series of conjectures the main
  of which is the log Minimal Model Program in dimension $n+1$.  These
  conjectures are theorems only when $n+1=3$, so only in the case of
  surfaces the results are not hypothetical, and this case was
  considered in detail in \cite{Alexeev_Mgn}.
  
  Where possible, we work in general context, over a fixed base
  scheme. The bulk of this material, however, applies only to the case
  of an algebraically closed field of characteristic $0$ because of
  the Minimal Model Program.
  
  We would like to point out that the general framework of what is
  described here has already been essentially understood in
  \cite{KollarShepherdBarron88}, \cite{Kollar_ProjCompleteModuli} and
  \cite{Alexeev_Mgn}. Many important ideas also come from
  \cite{Viehweg95}.
\end{say}

\medskip

\begin{say}
  We first remind what a stable $n$-pointed curve is.
\end{say}

\begin{defn}
  A stable $n$-pointed curve over an algebraically closed field is a
  collection $(C;P_1\dots P_m)$, where
  \begin{enumerate}
    \setcounter{enumi}{-1}
    \renewcommand{\theenumi}{(\arabic{enumi})}
  \item $C$ is a connected projective curve and $P_1\dots P_m$ are
    points on $C$.
  \item (condition on singularities) $C$ is reduced and has nodes
    only, and $P_1\dots P_m$ all lie in the nonsingular part.
  \item (numerical condition) for every smooth rational curve
    $E\subset C$, $E$ has at least 3 special points: one of $P_i$ or
    the nodes; and for every smooth elliptic curve or a rational curve
    with one node $E\subset C$, $E$ has at least 1 special point.
  \end{enumerate}
  
  A stable $n$-pointed curve over a scheme $S$ is a flat projective
  morphism $\pi:(C;P_1\dots P_m)\to S$, with $P_i\subset C$ closed
  subschemes and each $P_i\to S$ also flat, whose every geometric
  fiber is a stable $n$-pointed curve over a field $k=\bar k$.
\end{defn}

\begin{say}
  The moduli stack of stable $n$-pointed curves is proper, and it is
  coarsely represented by a projective scheme $M_{g,n}$, see
  \cite{Knudsen83}. 
\end{say}

\medskip

\begin{say}
  {\bf Question.} What is the analog of this in higher dimensions? 

  \medskip
  
  One definitely has to consider a collection consisting of a
  connected projective scheme $X$ plus $m$ closed subschemes. We have
  two basic choices: they could be points or divisors. Here, we choose
  divisors: $B_1\dots B_m$.

  The numerical condition (2) above can be reformulated by saying
  ``$K_C+\sum P_i$ is ample''. We can now directly transfer this to
  dimension $n$ if we understand what $K+B=K_X+\sum B_j$ is.

  Finally, the condition on the singularities. This is the trickiest
  of the three. The answer comes from the log Minimal Model Program
  theory: the singularities of $(X,B)$ have to be semi log canonical.
\end{say}

\medskip

\begin{say}
  We are now ready to introduce our main object.
\end{say}

\begin{defn}
  A stable pair over an algebraically closed field is a collection
  $(X;B_1\dots B_m)$, where
  \begin{enumerate}\setcounter{enumi}{-1}
    \renewcommand{\theenumi}{(\arabic{enumi})}
    \setcounter{enumi}{-1}
  \item $X$ is a connected projective not necessarily irreducible
    variety and $B_1\dots B_m$ are reduced divisors on $X$.
  \item (condition on singularities) the pair $(X,B)$ has semi log
    canonical singularities.
  \item (numerical condition) $K+B$ is ample.
  \end{enumerate}
  
  A stable pair over a scheme $S$ of level $N$ is a flat projective
  morphism $\pi:(X;B_1\dots B_m;\cL)\to S$, with $B_i\subset X$ closed
  subschemes, each $B_i\to S$ also flat and $\cL$ an invertible sheaf
  on $X$, whose every geometric fiber is a stable pair over a field
  $k=\bar k$ and such that the restriction of $\cL$ on each geometric
  fiber coincides with $\cO(N(K+B))$. We say that two pairs
  $(X_1,B_1;\cL_1)$ and $(X_2,B_2;\cL_2)$ are isomorphic if there
  exists an isomorphism of $(X_1,B_1)$ and $(X_2,B_2)$ over $S$ that
  induces a fiber-wise isomorphism of $\cL_1$ and $\cL_2$.
\end{defn}

\begin{conj}[Boundedness Conjecture]
  For every positive rational number $C$ there exist
  \begin{enumerate}
  \item a positive integer $N>0$ with the property that for every
    stable $n$-dimensional stable pair $(X,B)$ with $(K+B)^n=C$ the
    sheaf $\cO(N(K+B))$ is invertible.
  \item a scheme $S$ of finite type over the base scheme and a flat
    projective family $(X;B_1\dots B_m)$ whose geometric fibers
    include all stable $n$-dimensional pairs of level $N$ with
    $(K+B)^n=C$.
  \end{enumerate}
\end{conj}

\begin{say}
  This has been shown to be true only in dimension 2
  (\cite{Alexeev_Boundedness}) and trivially in dimension 1.
\end{say}

\begin{defn}
  We now fix a rational number $C$ and an integer $N$ as above and
  define the functor
  \begin{displaymath}
    \cM^N_C(S)=\left\{
      \begin{aligned}
        \text{stable $n$-dimensional pairs over }S \\
        \text{of level $N$ with } (K+B)^n=C 
      \end{aligned}
      \right\}/\simeq
  \end{displaymath}
  and the moduli stack by the same formula but without dividing by
  isomorphisms, and by giving $\cM^N_C(S)$ the groupoid structure in a
  natural way.
\end{defn}

\begin{say}
  There are other possible definitions for the moduli functor, see f.e
  \cite{Alexeev_Mgn}.
\end{say}

\medskip

\begin{say}
  At this point we can choose a certain scheme in a product of Hilbert
  schemes with the universal family that contains all interesting for
  us stable pairs.  The next step is to separate the stable pairs from
  wrong fibers, and for this we need to know that our functor is
  locally closed in the following sense.
  
  For every flat projective family $(X,B)\to S$ there exist
  locally closed subschemes $S_l\subset S$ with the following
  universal property:
  \begin{itemize}
  \item A morphism of schemes $T\to S$ factors through $\coprod
    S_l$ iff $(X,\cL)\underset{S}{\times}T \to T$ belongs to
    $\cM^N_C(T)$.
  \end{itemize}

  \medskip
  
  For our functor this property follows from the following conjecture
  of Shokurov (\cite{Shokurov91}).
\end{say}

\begin{conj}[Inversion of log Adjunction]
  \label{conj:inversion of log adjunction}
  Let $(X,B)\to S$ be a flat $1$-dimensional family. Assume that there
  exists an invertible sheaf $\cL$ on $X$ whose restriction on each
  fiber coincides with $\cO_X(N(K_X+B))$ as in definition
  \ref{defn:semi_log_canonical}. Then the $S_2$-fication of the pair
  $(X_0,B_0)$ has semi log canonical singularities iff the pair
  $(X,B+X_0)$ has semi log canonical singularities in a neighborhood
  of $X_0$. 
\end{conj}

\begin{rem}
  $X_0$ is $S_2$ iff $X$ is $S_3$. In many cases the varieties are
  Cohen-Macaulay, so taking the $S_2$-fication is unnecessary.
\end{rem}

\begin{saynum}\label{saynum:inversion_of_log_adjunction}
  One direction of this conjecture (going from the family to the
  central fiber) is easy and the proof for the case when $X$ is
  irreducible can be found in \cite[ch.17]{FAAT}. The general case can
  be easily deduced from that by taking the normalization.
  
  The same reference contains the proof of the opposite direction (the
  inversion) assuming the log Minimal Model Program in dimension
  $n+1$. It also contains several special cases where it can be proved
  without log MMP, using the Kawamata-Viehweg vanishing theorem only.
\end{saynum}

\medskip

\begin{saynum}
  The inversion of adjunction conjecture implies that if the sheaves
  $\cO_X(N(K_X+B))$ are locally free and can be put together in a flat
  family then the semi log canonical property is stable under
  generizations.
  
  Indeed, it follows from the definition that for a general fiber
  $X_t$ the pair $(X,B+X_0+X_t)$ is still semi log canonical. Then
  $(X_t,B_t)$ is semi log canonical by the easy direction of log
  adjunction.
  
  The question when exactly the sheaves $\cO_X(N(K_X+B))$ can be put
  together in a flat family is rather delicate. It follows from a
  technical result of Koll\'ar, see f.e. \cite{Alexeev_Mgn}.
\end{saynum}

\medskip

\begin{say}
  At this point we can pick a sub-family in our universal family that
  contains exactly our stable pairs. What remains is to take a
  quotient by the pre equivalence relation (or a groupoid) which is
  given by the action of the projective linear group $PGL$. This
  groupoid is easily seen to be flat. It also has a quasifinite
  stabilizer because stable pairs have finite automorphism groups by
  \cite{Iitaka82}. The next separateness property implies that the
  stabilizer is in fact finite. In this situation the quotient exists
  as a separated algebraic space. Nowadays, there are several
  convenient references for this statement, for example
  \cite{Kollar_QuotSpaces} and \cite{MoriKeel95}. As a result, one obtains a
  coarse moduli space $M_C^N$ as a separated algebraic space of finite
  type, and we are already working over an algebraically closed field
  $k$ of $\chr0$ since we used the log MMP.
\end{say}

\begin{thm}
  Let us assume the inversion of log adjunction conjecture. Let
  $(X',B')\to S\setminus 0$ be a $1$-dimensional family without the
  central fiber which is a stable pair over $S\setminus 0$. Then it
  can be completed to a stable pair over $S$ in no more than one way
  up to an isomorphism.
\end{thm}
\begin{proof}
  Let $(X,B)\to S$ be one such completion. By the inversion of log
  adjunction we know that $(X,B+X_0)$ is semi log canonical, possibly
  after shrinking $S$. Assume first that the scheme $X$ is
  irreducible, so that the singularities are in fact canonical. Then
  for any proper birational morphism from a normal variety $f:Y\to X$
  and for every positive integer $d$ we have by definition
  \begin{displaymath}
    f_*\cO_Y\big( dN(K_Y+f^{-1}B+f^{-1}X_0+\sum E_i) \big) 
    = \cO_X\big(dN(K_X+B+X_0)\big)
  \end{displaymath}
  Here the following three circumstances are important:
  \begin{enumerate}\renewcommand{\theenumi}{(\arabic{enumi})}
  \item $f^{-1}X_0+\sum E_i$ is in fact the central fiber of $Y$ with
    the reduced structure.
  \item the divisor $X_0$ is relatively trivial.
  \item the divisor $K+B$ is relatively ample, so that the family
    $(X,B)\to S$ can be computed as a $\Proj$ of a big graded ring of
    relative sections of $\cO(dN(K+B))$. 
  \end{enumerate}
  As a result of this, we obtain
  \begin{displaymath}
    X=\Proj_{d\ge0}\oplus\pi_*\cO_Y(dN(K_Y+f^{-1}B+Y_{0,red})),
  \end{displaymath}
  where $\pi$ denotes the morphism $Y\to S$.
  
  But this means that the family $(X,B)$ can be uniquely reconstructed
  from $(Y,f^{-1}B)$. Now, given two families $(X_1,B_1)$ and
  $(X_2,B_2)$. we can find a normal variety $Y$ which dominates both
  of them. By uniqueness, we have a canonical isomorphism
  $(X_1,B_1)\to(X_2,B_2)$.
  
  This completes the case when $X$ is irreducible. In general, the
  above argument shows the uniqueness of $(X^{\nu},B+\cond(\nu))$, and
  $(X,B)$ is uniquely recoverable from that.
\end{proof}

\begin{say}
  Next, we would like to prove that this algebraic space is in fact
  proper. For this, we have to check the corresponding property for
  our functor $\cM_C^N$.
\end{say}

\medskip

\begin{saynum}
  The pair $(X,B)$ above is the log canonical model of
  $(Y,f^{-1}B+Y_{0,red})$. So, the argument actually followed from the
  uniqueness of the log canonical model. Vice versa, assume that we
  have the log Minimal Model available. Start with arbitrary
  compactification $(X,B)\to S$ of a stable pair over $S\setminus0$.
  Take the normalization. For each irreducible component apply the
  Semistable Reduction Theorem (of course, $\chr 0$ is necessary for
  that) to obtain, after a finite ramified base change and resolution
  of singularities, a family with the reduced central fiber such that
  the irreducible components of the central fiber, $\cond(\nu)$,
  exceptional divisors of resolution and $B_j$ intersect
  transversally. Note that it is possible to choose the same base
  change that works for every irreducible component.  And then just
  find the log canonical model applying log MMP.
  
  In fact, we don't need all the results of log MMP but only the
  following conjecture and only in the $1$-dimensional semistable
  case. After that, glue the irreducible components back together.
  That will be the desired family over a finite ramified cover of $S$.
  This proves that our functor and the moduli space are proper.
\end{saynum}

\begin{conj}[Existence of log Canonical Model]
  Let $\pi:(Y,B)\to S$ be a projective morphism and assume that
  \begin{enumerate}
  \item the singularities of $(Y,B)$ are log canonical.
  \item restriction of $\cO_Y(N(K+B)$ on each generic fiber is big
    (contains an ample divisor).
  \end{enumerate}
  Then the ring of $\cO_S$-modules
  \begin{displaymath}
    \oplus_{d\ge0}\pi_*\cO_Y(dN(K_Y+B))
  \end{displaymath}
  is finitely generated.
\end{conj}

\begin{saynum}\label{saynum:ampleness_lemma}
  The last step is to show that the moduli space $M_C^N$ is
  projective. This follows by the Koll\'ar's Ampleness Lemma, see
  \cite{Kollar_ProjCompleteModuli}. The input data for this statement is
  \begin{enumerate}
  \item $M$ has to be a proper algebraic space of finite type over an
    algebraically closed field field $k$ of characteristic $0$.
  \item On a finite cover of $M$ there has to exist a projective
    polarized family $(X,B)$ whose every fiber has semi log canonical
    singularities (Koll\'ar considered the case $B=\emptyset$ but the
    generalization to the case of reduced $B$ is immediate). For
    example, this happens when $M$ is a coarse moduli space for some
    functor of polarized varieties, as in our case.
  \item The polarization has to be functorial, i.e. compatible with
    base changes. In our case, the polarization $\cO(N(K_{X/S}+B))$
    has this property.
  \end{enumerate}
\end{saynum}


\section{Examples of log canonical singularities}
\label{sec:Examples of log canonical singularities}

\begin{say}
  The following examples should be in any introductory article on log
  canonical singularities but surprisingly they aren't.
\end{say}

\begin{lem}\label{lem:torus_embeddings_are_slc}
  Let $X=T_Nemb(\Delta)$ be a torus embedding over a field $k$ defined
  by a rational partial polyhedral cone decomposition and $B=\sum B_j$
  be the sum of divisors corresponding to the $1$-dimensional faces of
  the fan $\Delta$. Then the pair $(X,B)$ has pre log canonical
  singularities, and $B$ (i.e. the pair $(B,0)$) has pre semi log
  canonical singularities.
\end{lem}
\begin{proof}
  The basic formula of the theory of torus embeddings for the
  canonical sheaf is
  \begin{displaymath}
    \omega_X(B)\simeq \cO_X
  \end{displaymath}
  Every torus embedding has a toric resolution of singularities
  $f:Y\to X$ such that $f^{-1}B\cup E_i$ has normal crossings, where
  $E_i$ are the exceptional divisors of $f$. Here $f^{-1}B\cup E_i$ is
  the union of divisors corresponding to $1$-dimensional faces of the
  fan of $Y$. Therefore,
  \begin{displaymath}
    f^*\cO(K_X+B)\simeq f^*\cO_X
    \simeq \cO_Y \simeq \cO(K_Y+f^{-1}B+\sum E_i)
  \end{displaymath}
  and the singularities of the pair $(X,B)$ are pre log
  canonical. 
  
  The normalization $B^{\nu}$ of $B$ is a disjoint union of torus
  embeddings, and $\cond(\nu)$ is again the union of divisors
  corresponding to the $1$-dimensional faces. This shows that $B$ has
  pre semi log canonical singularities.
\end{proof}

\begin{say}
  Therefore, every time when toric geometry is used, log canonical
  singularities show up. One of such situations is the following
  theorem of Mumford \cite[3.4,4.2]{Mumford_HirzebruchProportionality}.
\end{say}

\begin{thm}
  Let $\Gamma$ be a neat arithmetic group acting on a bounded
  symmetric complex domain $D$. Let $(D/\Gamma)^*$ be the Baily-Borel
  compactification of $D/\Gamma$ and $\overline{D/\Gamma}$ be any of
  the toroidal compactifications. Denote the boundaries of these
  compactifications by $\Delta^*$, $\overline{\Delta}$ respectively.
  Then
  \begin{align*}
    (D/\Gamma)^* &= 
    \Proj_{d\ge0} H^0\big( d(K_{(D/\Gamma)^*}+\Delta^*) \big) \\
    &= \Proj_{d\ge0} H^0\big( 
    d(K_{\overline{D/\Gamma}}+\overline{\Delta}) \big)
  \end{align*}
\end{thm}

\begin{cor}
  $((D/\Gamma)^*,\Delta^*)$ is the log canonical model of
  $(\overline{D/\Gamma},\overline{\Delta})$, and they both have log
  canonical singularities.
\end{cor}

\begin{say}
  The above formula in fact is one of the {\em definitions\/} of a log
  canonical model. We remind that a group $\Gamma$ is called neat if
  eigenvalues of each element of $\Gamma$ generate a torsion-free
  subgroup of $\bC^*$.  The quotient space $D/\Gamma$ by a neat group
  is nonsingular.
  
  What about the general case? It is easy, all one has to do is use
  the Hurwitz formula (cf. \cite[3.16]{Kollar_SingsPairs}).

  \medskip
  
  Every arithmetic group contains a neat subgroup
  $\Gamma_0\subset\Gamma$ of finite index. Let $D_j$ be the
  irreducible ramification divisors of
  $\overline{D/\Gamma_0}\to\overline{D/\Gamma}$ on
  $\overline{D/\Gamma}$ with ramification indices $n_j$. Then we
  immediately obtain the following
\end{say}

\begin{thm}
  \begin{align*}
    (D/\Gamma)^* &= 
    \Proj_{d\ge0} H^0\big( d(K_{(D/\Gamma)^*}+\Delta^*
    +\sum(1-1/n_j)D_j) \big) \\
    &= \Proj_{d\ge0} H^0\big( 
    d(K_{\overline{D/\Gamma}}+\overline{\Delta}
    +\sum(1-1/n_j)D^*_j) \big)
  \end{align*}
\end{thm}

\begin{cor}
  $((D/\Gamma)^*,\Delta^*+\sum(1-1/n_j)D^*_j)$ is the log canonical
  model of $(\overline{D/\Gamma},\overline{\Delta}+\sum(1-1/n_j)D_j)$,
  and they both have log canonical singularities.
\end{cor}

\begin{exmp}
  The compactification $\overline{A}_1=\bP^1$ of the moduli space
  $A_1$ of elliptic curves does not have log general type:
  \begin{displaymath}
    \deg(K_{\bP^1}+P_{\infty})= -2+1<0,  
  \end{displaymath}
  so it is not a log canonical model of anything.  However, the sum
  becomes positive when one adds the terms $(1-1/n_i)P_i$
  corresponding to the elliptic curves with automorphisms. This
  answers the footnote of Mumford appearing on the same page as
  theorem \cite[4.2]{Mumford_HirzebruchProportionality}.
\end{exmp}

\begin{say}
  Another situation is the stable quasiabelian varieties and pairs
  appearing as the limits of abelian varieties. We refer the reader to
  \cite{AlexeevNakamura96,Alexeev_CMAV} for their definition. The very
  construction for them is toric, so not surprisingly we have
\end{say}

\begin{lem}\label{lem:sqavs_are_slc}
  Let $P_0$ is a SQAV. Then $P_0$ has pre semi log canonical
  singularities.
\end{lem}
\begin{proof}
  By construction (\cite{AlexeevNakamura96}) there exists an \'etale
  map $\wP_0\to P_0$, and $\wP_0$ is a union of divisors in a torus
  embedding $\wP$ corresponding to the $1$-dimensional faces of the
  fan. The statement now follows from
  \ref{lem:torus_embeddings_are_slc}.
\end{proof}

\begin{say}
  $P_0$ in \cite{AlexeevNakamura96} appears as a central fiber of a
  one-dimensional degenerating normal family $P/S$ of abelian
  varieties. Over $\bC$, $P$ is a quotient of a torus embedding (which
  is locally of finite type) by a group $\bZ^g$ acting freely in the
  classic topology.
\end{say}

\begin{lem}
  The family $P$ itself has log canonical singularities.
\end{lem}
\begin{proof}
  Indeed, the general fiber of $P/S$ is smooth, so all the ``bad''
  discrepancies lie over the central fiber. By
  \ref{lem:torus_embeddings_are_slc} the pair $(P,P_0)$ has log
  canonical singularities, i.e. the corresponding discrepancies are
  $a_i\le1$. But the discrepancies of $(P,0)$ have to be less than
  $a_i$ by at least the multiplicities of $f^*P_0$ along the
  exceptional divisors. Since $P_0$ is Cartier, these multiplicities
  are at $\ge1$ and the discrepancies of $P$ are $\le0$.
\end{proof}

\begin{say}
  In the principally polarized case a SQAV by \cite{AlexeevNakamura96}
  comes with a natural theta divisor $\Theta$.
\end{say}

\begin{rem}
  An easy generalization of the last lemma is that a pair
  $(P_0,\varepsilon\Theta_0)$ has semi log canonical singularities for
  $\varepsilon\ll1$ in $\chr0$. For this one simply has to notice that
  $\Theta$ does not entirely contain any of the strata of $P_0$:
  \cite[3.28]{AlexeevNakamura96}.  A more interesting is the
  following.
\end{rem}

\begin{thm}
  \label{thm:sqaps_are_slc}
  A principally polarized stable quasiabelian pair $(P_0,\Theta_0)$
  over $\bC$ has semi log canonical singularities.
\end{thm}
\begin{proof}
  For the abelian varieties this result is a theorem of Koll\'ar
  \cite{Kollar_ShafMapsnPlurigenera}. The present proof is the
  adaptation of the proof of that theorem to our situation.
  
  By \cite{AlexeevNakamura96} every stable quasiabelian pair appears
  as the central fiber in a $1$-dimensional family $\pi:(P,\Theta)\to
  S=D_{\varepsilon}$ with abelian general fiber over a small disk.  We
  denote by $I$ the ideal defining $0\in D_{\varepsilon}$.
  
  If we prove that the pair $(P,\Theta+P_0)$ has log canonical
  singularities then we would be done by the easy direction of the
  ``inversion of log adjunction theorem'' (see \cite[ch.17]{FAAT} or
  \ref{conj:inversion of log adjunction}).
  
  The locus $Z$ of non-log canonical singularities of $(P,\Theta+P_0)$
  coincides with the locus of non-log terminal singularities of the
  pair $(P,(1-\varepsilon)(\Theta+P_0))$ for $0<\varepsilon\ll1$. We
  will apply the Kawamata-Viehweg vanishing theorem in the following
  Nadel's form (see f.e. \cite[2.16]{Kollar_SingsPairs}):

  \begin{thm}
    Let $X$ be a normal and proper variety and $N$ a line bundle on
    $X$. Assume that $N\equiv K_X+\Delta+M$, where $M$ is nef and big
    $\bQ$-Cartier divisor and $\Delta$ effective $\bQ$-Cartier divisor
    with coefficients $<1$. Then there is an ideal sheaf
    $J\subset\cO_X$ such that
    \begin{displaymath}
      \Supp(\cO_X/J)=\{ x\in X \,|\, (X,\Delta)
      \text{ is not log terminal at } x \}
    \end{displaymath}
  \end{thm}

  We will apply this theorem in the relative situation to the proper
  morphism $\pi:P\to S$. We have 
  \begin{displaymath}
    K_P+\Theta+P_0=K_P+(1-\varepsilon)(\Theta+P_0)+
    \varepsilon(\Theta+P_0),
  \end{displaymath}
  $K_P,P_0$ are relatively trivial and $\Theta$ is relatively ample.
  Therefore, by the above there exists an ideal $J\subset\cO_P$
  supporting the locus $Z$ where the pair $(P,\Theta+P_0)$ is not log
  canonical, and $R^1\pi_*J(\Theta+P_0)=0$. Therefore, the following
  map is surjective
  \begin{displaymath}
    \pi_*\cO_P(K_P+\Theta+P_0)
    \underset{\phi}{\to} \pi_*\cO_Z(K_P+\Theta+P_0)
  \end{displaymath}
  
  Since in the nonsingular case the statement holds by Kollar's
  theorem, after shrinking $S$ the support of $Z_{red}$ will be
  contained in the central fiber. Therefore $Z$ is a closed
  complex-analytic subspace of $P_n=P\underset{R}{\times}R/I^n$ for
  some $n\ge0$, where $R$ is the ring of germs of analytic functions
  at $0$.  Moreover, $Z$ is a closed subspace of the theta
  divisor $\Theta$.  Indeed, as in the proof of lemma
  \ref{lem:sqavs_are_slc}, theorem \ref{lem:torus_embeddings_are_slc}
  the pair $(P,P_0)$ is log canonical, therefore the pair
  $(P,(1-\varepsilon)P_0)$ is log terminal.
  
  By theorem 4.6 of \cite{AlexeevNakamura96} we have
  $H^0(\cO_{P_0}(\Theta_0))=1$ and $H^i(\cO_{P_0}(\Theta_0))=0$ for
  $i>0$. This implies $R^i\pi_*\cO_P(\Theta)=0$ for $i>0$ and
  $\pi_*\cO_P(\Theta)=\cO_S$.
  
  We have $\cO(K_P)\simeq\cO(P_0)\simeq\cO$, so
  $\cO(K_P+\Theta+P_0)\simeq \cO(\Theta)$.  Since $Z$ is a closed
  subspace of $\Theta$, $\phi$ has to be the zero map.  On the other
  hand, $\pi_*\cO_Z(\Theta)\ne0$ for any proper subspace $Z\subset
  P_n$.  In the nonsingular case this is concluded by the
  semi-continuity argument and the fact that the abelian variety acts
  transitively by translations. In our situation, there is the action
  of a semiabelian group $G/S$, and although it is not transitive,
  still the intersection of translations $g(\Theta)$ by sections $g\in
  G$ is empty: \cite[3.28]{AlexeevNakamura96}, and this implies
  $\pi_*\cO_Z(\Theta)\ne0$.
\end{proof}

\begin{cor}
  Over $\bC$, if the compactification of the moduli space $A_g$ by the
  pairs $(P_0,\Theta_0)$ exists, it is projective.
\end{cor}
\begin{proof}
  This follows by applying the Ampleness Lemma of Koll\'ar, cf.
  \ref{saynum:ampleness_lemma}.
\end{proof}


\ifx\undefined\bysame
\newcommand{\bysame}{\leavevmode\hbox to3em{\hrulefill}\,}
\fi


\begin{thebibliography}{{Kol}92}

\bibitem[Ale94]{Alexeev_Boundedness}
V.~Alexeev, {\em Boundedness and ${K}^2$ for log surfaces}, Int. J. Math. {\bf
  5} (1994), no.~6, 779--810.

\bibitem[Ale96a]{Alexeev_CMAV}
V.~Alexeev, {\em Compact moduli of (co)abelian varieties}, Preprint (1996).

\bibitem[Ale96b]{Alexeev_Mgn}
V.~Alexeev, {\em Moduli spaces ${M}_{g,n}$ for surfaces}, (M.~Andreatta and
  T.~Peternell, eds.), Walter de Gruyter, 1996, Proceedings of the
  international conference held in Trento, Italy, June 15-24, 1994, pp.~1--22.

\bibitem[AN96]{AlexeevNakamura96}
V.~Alexeev and I.~Nakamura, {\em On {M}umford's construction of degenerating
  abelian varieties}, Preprint (1996).

\bibitem[Har77]{Hartshorne77}
R.~Hartshorne, {\em Algebraic geometry}, Graduate Texts in Mathematics,
  vol.~52, Springer-Verlag, 1977.

\bibitem[Iit82]{Iitaka82}
S.~Iitaka, {\em Algebraic geometry. {A}n introduction to birational geometry of
  algebraic varieties}, Graduate Texts in Mathematics, vol.~76, Springer-Verlag
  New York, Inc., 1982.

\bibitem[KM95]{MoriKeel95}
S.~Keel and S.~Mori, {\em Quotients by groupoids}, Preprint (1995).

\bibitem[Knu83]{Knudsen83}
F.F. Knudsen, {\em The projectivity of the moduli space of stable curves, {II}:
  the stacks ${M}_{g,n}$}, Math. Scand. {\bf 52} (1983), 161--199.

\bibitem[Kol]{Kollar_SingsPairs}
J.~Koll\'ar, {\em Singularities of pairs}, Preprint.

\bibitem[Kol90]{Kollar_ProjCompleteModuli}
J.~Koll\'ar, {\em Projectivity of complete moduli}, J. Diff. Geom. {\bf 32}
  (1990), 235--268.

\bibitem[{Kol}92]{FAAT}
J.~{Koll\'ar et al.}, {\em Flips and abundance for algebraic threefolds},
  Ast{\'e}risque {\bf 211} (1992), 1--258.

\bibitem[Kol93]{Kollar_ShafMapsnPlurigenera}
J.~Koll{\'a}r, {\em Shafarevich maps and plurigenera of algebraic varieties},
  Invent. Math. {\bf 113} (1993), no.~1, 177--215.

\bibitem[Kol95]{Kollar_QuotSpaces}
J.~Koll{\'a}r, {\em Quotient spaces modulo algebraic groups}, Preprint (1995).

\bibitem[KSB88]{KollarShepherdBarron88}
J.~Koll\'ar and N.~Shepherd-Barron, {\em Threefolds and deformations of surface
  singularities}, Invent. Math. {\bf 91} (1988), 299--338.

\bibitem[Mum77]{Mumford_HirzebruchProportionality}
D.~Mumford, {\em Hirzebruch's proportionality theorem in the non-compact case},
  Invent. Math. {\bf 42} (1977), 239--272.

\bibitem[Sho92]{Shokurov91}
V.V. Shokurov, {\em 3-fold log flips}, Math. USSR--Izv. {\bf 56} (1992),
  105--203.

\bibitem[Vie95]{Viehweg95}
E.~Viehweg, {\em Quasi-projective moduli for polarized manifolds},
  Springer-Verlag, 1995.

\end{thebibliography}
\end{document}